\documentclass{article}

\usepackage{zharticle}
\usepackage{graphicx,amsmath,amsfonts,verbatim}
\usepackage{tikz}
\numberwithin{equation}{section}
\allowdisplaybreaks%
\usepackage{mathtools}
\usepackage{enumitem}

\newcounter{algorithm}[section]
\renewcommand{\thealgorithm}{\thesection.\arabic{algorithm}}
\newenvironment{algorithm}[1]{%
    \refstepcounter{algorithm}%
    \begin{list}{}{%
		\setlength{\rightmargin}{0\linewidth}%
		\setlength{\leftmargin}{.05\linewidth}%
    }%
		\small%
		\item[]{\setlength{\parskip}{0ex}\hrulefill\par\nopagebreak%
        {\sffamily\bfseries Algorithm~\thealgorithm}\ {\scshape#1}}%
        \begin{itemize}[label={}, leftmargin=0pt]%
}{%
    \end{itemize}%
    \vspace{-\parskip}
    {\setlength{\parskip}{-1ex}\nopagebreak\par\hrulefill}%
    \end{list}
}

\newcommand{\ones}{{\bf 1{}}}  
\newcommand{\reals}{{\bf R{}}}  
\newcommand{\symms}{{\bf S{}}}  
\newcommand{\prob}{\mathop{\bf prob{}}}  
\newcommand{\dom}{\mathop{\bf dom}}  
\newcommand{\tr}{\mathop{\bf tr}}  
\newcommand{\argmax}{\mathop{\rm argmax}}  
\newcommand{\argmin}{\mathop{\rm argmin}}  
\newcommand{\norm}[1]{{\lVert #1 \rVert}}  
\newcommand{\cone}{{\cal K{}}}  
\newcommand{\dset}{{\cal D{}}}  
\newcommand{\cset}{{\cal C{}}}  
\newcommand{\xset}{{\cal X{}}}  
\newcommand{\yset}{{\cal Y{}}}  
\newcommand{\Dkl}{\mathop{D_{\rm kl}}}  

\newcommand{\cf}{{\rmfamily\itshape cf.}}
\newcommand{\eg}{{\rmfamily\itshape e.g.}}
\newcommand{\ie}{{\rmfamily\itshape i.e.}}

\newcommand{\etal}{{\rmfamily\itshape et al.}}

\usepackage{lipsum}

\usepackage[colorlinks]{hyperref}

\usepackage{listings}
\usepackage{xcolor}

\definecolor{codegreen}{rgb}{0,0.6,0}
\definecolor{codegray}{rgb}{0.5,0.5,0.5}
\definecolor{codepurple}{rgb}{0.58,0,0.82}
\definecolor{backcolour}{rgb}{0.95,0.95,0.98}

\lstdefinestyle{zhstyle}{
    backgroundcolor=\color{backcolour},
    commentstyle=\color{codegreen},
    keywordstyle=\color{magenta},
    numberstyle=\tiny\color{codegray},
    stringstyle=\color{codepurple},
    basicstyle=\ttfamily\small,
    breakatwhitespace=false,
    breaklines=true,
    captionpos=b,
    keepspaces=true,
    numbers=left,
    numbersep=5pt,
    showspaces=false,
    showstringspaces=false,
    showtabs=false,
    tabsize=2,
}
\lstdefinelanguage{zhpython}{%
    sensitive=true,%
    morekeywords={and, as, assert, break, class, continue, def, del, elif,%
        else, except, exec, finally, for, from, global, if, import, in, is,%
        lambda, not, or, pass, print, raise, return, try, while, with, yield},%
    morecomment=[l]\#,%
    morestring=[s]{'''}{'''},%
    morestring=[s]{"""}{"""},%
    morestring=[b]',%
    morestring=[b]"%
}

\title{Disciplined Biconvex Programming}
\author[1,2]{Hao Zhu}
\author[1,2]{Joschka Boedecker}
\affil[1]{IMBIT//BrainLinks-BrainTools}
\affil[2]{Department of Computer Science, University of Freiburg}

\begin{document}
\maketitle

\begin{abstract}
    We introduce \emph{disciplined biconvex programming} (DBCP), a modeling framework for specifying and solving biconvex optimization problems.
    Biconvex optimization problems arise in various applications, including machine learning, signal processing, computational science, and control.
    Solving a biconvex optimization problem in practice usually involves heuristic methods based on alternate convex search (ACS), which iteratively optimizes over one block of variables while keeping the other fixed, so that the resulting subproblems are convex and can be efficiently solved.
    However, designing and implementing an ACS solver for a specific biconvex optimization problem usually requires significant effort from the user, which can be tedious and error-prone.
    DBCP extends the principles of disciplined convex programming to biconvex problems, allowing users to specify biconvex optimization problems in a natural way based on a small number of syntax rules.
    The resulting problem can then be automatically split and transformed into convex subproblems, for which a customized ACS solver is then generated and applied.
    DBCP allows users to quickly experiment with different biconvex problem formulations, without expertise in convex optimization.
    We implement DBCP in the open-source Python package \texttt{dbcp}, as an extension to the well known domain specific language CVXPY for convex optimization.
\end{abstract}

\vspace*{\fill}
\begin{center}\sffamily
    This manuscript has been accepted for publication in \emph{Journal of Optimization Theory and Applications}.
    This preprint version should include all the content (subject to minor changes) in the final published version, but may differ in formatting and copyediting.
\end{center}
\vspace*{\fill}

\newpage
\tableofcontents
\newpage

\section{Introduction}
We consider biconvex optimization problems, which involve optimizing a biconvex objective function subject to a biconvex set constraint.
Informally, by saying that a set (or function) is biconvex, we mean that it is convex in each of two blocks of variables when the other block is fixed.
Biconvex optimization problems arise in various applications, including machine learning~\cite{paatero1994positive,udell2016generalized,zhu2025multi}, signal and image processing~\cite{levin2009understanding,chaudhuri2014blind}, recommender systems~\cite{chu2009personalized,johnson2014logistic}, control~\cite{safonov1994control,vanantwerp2000tutorial,javanmardi2020bilinear}, and brain computer interfaces~\cite{dyrholm2007bilinear,shi2014sparse}.

Unlike convex optimization problems, biconvex optimization problems are in general nonconvex and can be very hard to solve.
In fact, many of them have been shown to be NP-hard (see \cite{wen2012solving}, \cite{toker1995np}, and \cite{zhu2025multi} for some examples).
Nevertheless, various heuristic methods based on \emph{alternate convex search} (ACS)~\cite{de1994block} exist for finding good solutions to many biconvex optimization problems in practice.
The basic idea of ACS-type methods is to iteratively optimize over one block of variables while keeping the other block fixed, such that each subproblem is convex and can be efficiently solved.
Although in the most general case, theoretical convergence of ACS-type methods is only guaranteed to stationary points, these methods turn out to work quite well in practice and hence have become the most popular choice for finding a satisfactory solution to biconvex optimization problems.

Designing and implementing an ACS solver for a specific biconvex optimization problem usually requires significant effort from the user, including properly partitioning the biconvex problem into convex subproblems, implementing robust and efficient solvers for the subproblems, and designing appropriate stopping criteria.
Moreover, maintaining numerical stability and ensuring convergence of the ACS procedure typically necessitate integrating appropriate regularizations or modifications into the subproblems to ensure properties such as strong convexity and feasibility.
With the help of \emph{domain specific languages} (DSLs) for convex optimization, such as CVXPY~\cite{diamond2016cvxpy,agrawal2018rewriting}, specifying and solving the (modified) convex subproblems is largely simplified and automated, while the other steps still require the user to have expert knowledge about convex analysis and programming, and can be tedious and error-prone.

In this paper, we propose a modeling framework for biconvex optimization problems, named \emph{disciplined biconvex programming} (DBCP), which extends the ideas of \emph{disciplined convex programming} (DCP)~\cite{grant2006disciplined} to biconvex problems.
DBCP for biconvex optimization is analogous to DCP for convex optimization, which allows users to specify biconvex optimization problems in a natural way based on a small number of syntax rules.
When an optimization problem description complies with the DBCP syntax rules, it is guaranteed to be a valid biconvex problem, and more importantly, it can be automatically split and transformed into convex (specifically, DCP-compliant) subproblems.
A customized ACS solver with procedures for handling initialization and improving numerical stability is then generated and applied.
Like DCP does for convex optimization, DBCP makes it easy to specify and solve biconvex problems, allowing users to quickly prototype and experiment with different biconvex problem formulations, without expertise in (or even knowledge of) convex analysis and solution methods for biconvex problems.
In fact, most users might be unaware of the variable partition, transformation, and solution processes, through which a solution to the biconvex problem is found.
We implement DBCP in a Python package, named \texttt{dbcp}, as an extension of CVXPY, which is fully open-source and available at
\begin{quote}
    \url{https://github.com/dxogrp/dbcp}.
\end{quote}

\subsection{Previous and related work}
\paragraph{Biconvex analysis.}
The first recognition of biconvex structure in the context of mathematical programming can be traced back to Falk and Soland~\cite{falk1969algorithm} in the 1960s.
Only a few papers exist in the literature where biconvex sets are investigated.
Most of the corresponding theoretical results can be found in the papers of Aumann and Hart~\cite{aumann1986bi} and Goh~\etal~\cite{goh1994biaffine}.
Biconvex functions appear regularly in practice and hence have been discussed widely in the literature.
Properties of biconvex functions that are most relevant to optimization problems can be found in the works of Goh~\etal~\cite{goh1994biaffine} and Gorski~\etal~\cite{gorski2007biconvex}.
In addition, biconvex functions play an important role in martingale theory, and, in particular, in the analysis of Banach and Hilbert spaces~\cite{burkholder1981geometrical,burkholder1986martingales,aumann1986bi,lee1993burkholder}.
Biconvex functions can also be used to derive results on the robust stability of control systems in practical control engineering~\cite{geng2000robusta,geng2000robustb}.
Thibault~\cite{thibault1984continuity}, Jouak and Thibault~\cite{jouak1985directional}, and Borwein~\cite{borwein1986partially} analyzed the continuity and differentiability of measurable biconvex operators in topological vector spaces.
Al-Khayyal and Falk~\cite{al1983jointly} published results on the maximization of biconvex functions.
In the most general case, very little can be said about the global or even local optimality properties of biconvex optimization problems.
Nevertheless, some useful properties of partial optimality conditions (\ie, conditions for stationary points) of biconvex problems are discussed in~\cite{wendell1976minimization,gorski2007biconvex}.

\paragraph{Solution methods for biconvex problems.}
Current solution methods for biconvex problems can be roughly categorized into two classes: heuristic methods and global optimization methods.
Heuristic methods for biconvex problems aim at finding stationary points of the biconvex objective function, and are mostly based on the idea of alternately optimizing two convex subproblems, \ie, the ACS-type methods.
ACS methods are a special case of \emph{block relaxation methods}~\cite{warga1963minimizing,powell1973search,de1994block} where the variables are divided into disjoint blocks, and in each step, only the variables of an active block are optimized while those of the other blocks are fixed.
In particular, for ACS, only two blocks of variables defined by the convex subproblems are activated in cycles.
Since the resulting subproblems are convex, efficient convex minimization methods can be used to solve these subproblems.
A survey on ACS methods for biconvex optimization problems can be found in~\cite{wendell1976minimization}.
Gorski~\etal~\cite{gorski2007biconvex} showed that under weak assumptions, all solution points generated by ACS form a compact connected set and that each of these points is a stationary point of the objective function.
However, no better convergence results regarding local or global optimality properties can be obtained in general.
Regarding the global optimization methods, Floudas and Visweswaran~\cite{floudas1990global} adapted the idea of \emph{branch-and-bound}~\cite{lawler1966branch,boyd2007branch} to solve biconvex problems to global optimality.
Detailed mathematical background and an outline of this algorithm are given in~\cite{floudas1990global,floudas1993primal,floudas1995nonlinear,floudas2000deterministic}, and some basic convergence properties of this method are discussed in~\cite{gorski2007biconvex}.

\paragraph{Biconvex optimization applications.}
The practical usefulness of biconvex programming in many applications is also increasingly well known.
In machine learning, nonnegative matrix factorization~\cite{paatero1994positive,lee1999learning,udell2016generalized}, dictionary learning~\cite{aharon2006k}, and bilinear regression~\cite{hu2008collaborative,chu2009personalized,johnson2014logistic} are all well known biconvex optimization problems, which have been widely applied in representation learning, signal processing, and recommender systems.
The very popular $k$-means clustering problem can also be formulated as a biconvex program~\cite{lloyd1982least,so2022convergence,zhu2025multi}.
Blind deconvolution problems~\cite{levin2009understanding,chaudhuri2014blind} are another important class of biconvex optimization problems for image processing and communication.
Other studies such as Mishra~\etal~\cite{mishra2017sequential} and Fosson~\cite{fosson2018biconvex} explored biconvex structures in sparse learning and reweighted regression.
Biconvex problems also appear when dealing with bilinear matrix inequalities for bilinear control system synthesis~\cite{safonov1994control,vanantwerp2000tutorial,javanmardi2020bilinear}.
In economics, Schaller and Boyd~\cite{schaller2026estimating} formulate the regularized maximum likelihood estimation of a price elasticity matrix with a low-rank-plus-diagonal structure as a biconvex problem in the two low rank factors.
Finally, Dyrholm~\etal~\cite{dyrholm2007bilinear} and Shi~\etal~\cite{shi2014sparse} applied bilinear discriminant component analysis techniques to brain computer interfaces.
We discuss some of these applications in more detail in \S\ref{sec:example}.

\paragraph{Disciplined convex programming.}
DCP~\cite{grant2006disciplined} is a modeling framework for convex optimization problems.
If a mathematical optimization problem is written following the DCP syntax ruleset, it is guaranteed to be a valid convex program, and can then be automatically transformed into an appropriate canonical form~\cite{agrawal2018rewriting} for a numerical convex optimization solver, such as OSQP~\cite{stellato2020osqp}, ECOS~\cite{domahidi2013ecos}, SCS~\cite{o2016conic}, and Clarabel~\cite{goulart2024clarabel}.
Many DSLs for convex optimization, such as YALMIP~\cite{lofberg2004yalmip}, CVX~\cite{grant2014cvx}, Convex.jl~\cite{udell2014convex}, CVXPY~\cite{diamond2016cvxpy}, and CVXR~\cite{fu2020cvxr}, are based on the DCP ruleset.
Over the years, DCP has been extended to handle stochastic~\cite{ali2015disciplined}, convex-concave~\cite{shen2016disciplined}, geometric~\cite{agrawal2019disciplined}, quasiconvex~\cite{agrawal2020disciplined}, and saddle~\cite{schiele2024disciplined} problems.
The most closely related extension of DCP to our work is \emph{disciplined multi-convex programming} (DMCP)~\cite{shen2017disciplined}, which provides a DSL for specifying and solving multi-convex optimization problems via block relaxation methods.
However, it seems to be no longer actively supported or maintained.
In this paper, we work with CVXPY and focus on biconvex problems only, since most multi-convex problems in practice can be formulated as biconvex problems by grouping variables appropriately.
\emph{Disciplined nonlinear programming} (DNLP)~\cite{cederberg2026disciplined} is a recent extension of CVXPY for specifying and solving general nonlinear programs, including some biconvex problems.
For more detailed comparisons of the proposed DBCP framework with DMCP and DNLP, see \S\ref{sec:comp_dmcp} and \S\ref{sec:comp_dnlp}, respectively.

\subsection{Outline}
This paper is \emph{not} about establishing new theoretical results or computational methods for biconvex optimization problems.
Instead, we collect well known ideas and assemble them into a DSL for specifying and solving biconvex optimization problems in a natural human-readable way which is close to the mathematical formulation.
These ideas are implemented in the open-source Python package \texttt{dbcp}, which extends CVXPY to support biconvex programming.

The rest of this paper is organized as follows.
In \S\ref{sec:bcp}, we provide formal definitions for biconvex sets, functions, and optimization problems, including some specific examples of these objects, and briefly review their basic properties.
In \S\ref{sec:bcp_solution}, we introduce the ACS-based method for solving biconvex optimization problems, and discuss some practical augmentations for handling initialization and improving numerical stability, which are all integrated into the \texttt{dbcp} package.
Then in \S\ref{sec:dbcp}, we introduce the DBCP biconvex syntax ruleset for specifying biconvex optimization problems, as an extension of the DCP ruleset.
The implementation and basic usage of the provided Python package \texttt{dbcp}, where the functions and features mentioned above are implemented, are discussed in \S\ref{sec:impl}.
Finally, in \S\ref{sec:example}, we present some specific numerical examples for specifying and solving biconvex optimization problems that frequently appear in practice using the \texttt{dbcp} package, to show the simplicity and effectiveness of our framework.
These examples along with the corresponding code snippets also serve as a preliminary tutorial for new users to get started with \texttt{dbcp}.

\section{Biconvex programming}\label{sec:bcp}
In this section, we provide formal definitions for biconvex sets, biconvex functions, and biconvex optimization problems to establish our notation, and we briefly review some of the basic properties of these objects.
To start with, let $\xset \subseteq \reals^n$ and $\yset \subseteq \reals^k$ be two nonempty closed convex sets.
We should also note that throughout this paper, we use the convention that functions are implicitly extended-valued, \ie, by writing $f \colon A \to \reals$, we really mean $f \colon A \to \reals \cup \{\infty\}$.
In other words, a function defined on $A$ can take the value $\infty$ for some points in $A$~\cite{rockafellar1970convex,boyd2004convex}.

\subsection{Biconvex sets}\label{sec:bcvxset}
A set $B \subseteq \xset \times \yset \subseteq \reals^{n+k}$ is a \emph{biconvex set} if for every fixed $y \in \yset$, the set
\begin{equation*}
    B_y = \{x \in \xset \mid (x,y) \in B\} \subseteq \reals^n
\end{equation*}
is convex, and for every fixed $x \in \xset$, the set
\begin{equation*}
    B_x = \{y \in \yset \mid (x,y) \in B\} \subseteq \reals^k
\end{equation*}
is convex.

\begin{figure}[t]
    \centering
    \includegraphics[width=0.9\textwidth]{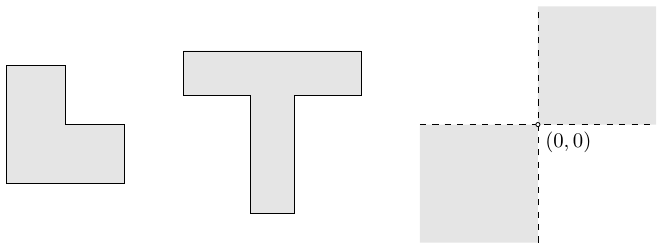}
    \caption{Examples of biconvex sets.}\label{fig:bcvxset}
\end{figure}

Obviously, a biconvex set is not necessarily convex in general.
Figure~\ref{fig:bcvxset} shows some examples of biconvex sets which are not convex.
In particular, a biconvex set does not even have to be connected, as shown by the example
\begin{equation*}
    B = \{(x, y) \in \reals^2 \mid x, y < 0\} \cup \{(x, y) \in \reals^2 \mid x, y > 0\},
\end{equation*}
whose picture is shown on the right of figure~\ref{fig:bcvxset}.

We should note that whenever several biconvex sets, biconvex functions, or biaffine mappings are considered together in the properties below, they are assumed to share the same fixed partition of variables into the $x$- and $y$-blocks.
An important property of biconvex sets is that, as with convex sets, the intersection of an arbitrary collection of biconvex sets is still a biconvex set~\cite{aumann1986bi,gorski2007biconvex}.

\subsection{Biconvex functions}\label{sec:bcvx_fn}
\begin{figure}
    \centering
    \includegraphics[width=0.6\textwidth]{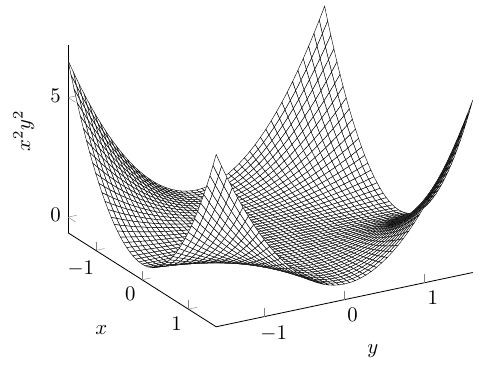}
    \caption{Graph of the biconvex function $f(x, y) = x^2 y^2$.}\label{fig:bcvx_fn}
\end{figure}

A function $f \colon \xset \times \yset \to \reals$ is a \emph{biconvex function} if its effective domain
\begin{equation*}
    \dom f = \{(x,y) \in \xset \times \yset \mid f(x,y) < \infty\}
\end{equation*}
is a biconvex set, and for every fixed $y \in \yset$, the function
\begin{equation*}
    f_y \colon \xset \to \reals,\quad x \mapsto f(x, y),
\end{equation*}
is convex in $x$, and for every fixed $x \in \xset$, the function
\begin{equation*}
    f_x \colon \yset \to \reals,\quad y \mapsto f(x, y),
\end{equation*}
is convex in $y$.
In other words, a biconvex function is one that is convex in each of two blocks of variables when the other block is fixed.
We can also define \emph{biconcave}, \emph{biaffine}, and \emph{bilinear} functions similarly, by replacing the property of being convex for $f_y$ and $f_x$ by the property of being concave, affine, or linear, respectively.
Figure~\ref{fig:bcvx_fn} shows an example of a biconvex function given by $f(x, y) = x^2 y^2$.

Now we list some basic properties of biconvex functions.
We do not provide formal proofs here;
interested readers may refer to \cite{goh1994biaffine} and \cite{gorski2007biconvex} for the corresponding proofs and discussion of the following statements (which are more or less obvious with the basic convex analysis toolbox~\cite{rockafellar1970convex}).
Similar to convex functions, the sublevel sets of a biconvex function $f \colon \xset \times \yset \to \reals$, given by
\begin{equation*}
    C_\alpha = \{(x, y) \in \dom f \mid f(x, y) \leq \alpha\},
\end{equation*}
are biconvex sets, for any $\alpha \in \reals$.
Also, many arithmetic properties that are valid for convex functions can be transferred to the biconvex case:
\begin{itemize}
    \item \emph{Nonnegative weighted sums.}
        Evidently, if $f$ is a biconvex function and $w \geq 0$, then the function $w f$ is biconvex.
        If $f_1$ and $f_2$ are both biconvex functions, then so is their sum $f_1 + f_2$.
        In the general case, if $f_1, \ldots, f_m$ are biconvex functions, and $w_1, \ldots, w_m \geq 0$, then the function
        \begin{equation*}
            f = w_1 f_1 + \cdots + w_m f_m
        \end{equation*}
        is also biconvex.
    \item \emph{Composition with a biaffine mapping.}
        If $h \colon \reals^m \to \reals$ is convex, and $g \colon \xset \times \yset \to \reals^m$ is biaffine, then the function $f \colon \xset \times \yset \to \reals$, given by $f(x, y) = h(g(x, y))$, is biconvex.
        In particular, if $h \colon \reals \to \reals$ is convex, and $A \in \reals^{m \times n}$, $C \in \reals^{m \times k}$, $b, d \in \reals^m$, then the function $f \colon \xset \times \yset \to \reals$, given by
        \begin{equation*}
            f(x, y) = h({(Ax + b)}^T(Cy + d))
        \end{equation*}
        is biconvex.
    \item \emph{Pointwise maximum and supremum.}
        If $f_1, \ldots, f_m$ are biconvex functions, then the function
        \begin{equation*}
            f(x, y) = \max\{f_1(x, y), \ldots, f_m(x, y)\}
        \end{equation*}
        is also biconvex.
        More generally, if ${\{f_i\}}_{i \in I}$ is a family of biconvex functions indexed by a set $I$, then the function
        \begin{equation*}
            f(x, y) = \sup_{i \in I} f_i(x, y)
        \end{equation*}
        is biconvex.
    \item \emph{Composition.}
        If $h \colon \reals \to \reals$ is convex and nondecreasing, and $g \colon \xset \times \yset \to \reals$ is biconvex, then the function $f \colon \xset \times \yset \to \reals$, given by $f(x, y) = h(g(x, y))$, is biconvex.
        This property can be easily extended to multivariate functions $h \colon \reals^m \to \reals$ that are convex and nondecreasing in each argument.
        Interested readers may refer to \cite[\S3.2.4]{boyd2004convex} for more details.
\end{itemize}
We should note that the property of \emph{composition with a biaffine mapping} above is a special case of the more general \emph{composition} property, but is of more practical importance, since most biconvexity structures in applications arise from composing convex functions with biaffine mappings.
(See \S\ref{sec:example} for some examples.)

\subsection{Biconvex optimization problems}\label{sec:bcvxprob}
A \emph{biconvex optimization problem} is one of the form
\begin{equation}\label{prob:bcp}
    \begin{array}{ll}
        \mbox{minimize} & f_0(x, y)\\
        \mbox{subject to} & f_i(x, y) \leq 0,\quad i = 1, \ldots, m\\
        & h_i(x, y) = 0,\quad i = 1, \ldots, p,
    \end{array}
\end{equation}
where $x \in \xset$, $y \in \yset$ are the optimization variables.
The functions $f_i \colon \xset \times \yset \to \reals$, $i = 0, \ldots, m$, are biconvex, and $h_i \colon \xset \times \yset \to \reals$, $i = 1, \ldots, p$, are biaffine.
Here (and also in \S\ref{sec:bcp_solution}), when discussing theoretical properties, we assume that the problem variables are partitioned into two disjoint and exhaustive blocks, $x$ and $y$, for notational convenience and in accordance with the convention in the literature.
However, in practice, the two groups passed to our package interface must remain disjoint but need not be exhaustive; see \S\ref{sec:impl} for more details.
Roughly speaking, the problem (\ref{prob:bcp}) can be interpreted as minimizing a biconvex objective function over a biconvex feasible set defined by biconvex inequality constraints and biaffine equality constraints, since each of these constraints defines a biconvex set, and the intersection of biconvex sets is still a biconvex set (\cf~\S\ref{sec:bcvxset}).
The \emph{feasible set} of (\ref{prob:bcp}) is defined as
\begin{equation*}
    \dset = \left\{(x, y) \in \xset \times \yset\ \middle|\ \begin{array}{l}
        f_0(x, y) < \infty\\
        f_i(x, y) \leq 0,\quad i = 1, \ldots, m\\
        h_i(x, y) = 0,\quad i = 1, \ldots, p
    \end{array}\right\}.
\end{equation*}

Different from convex optimization problems, in the most general case, very little can be said about the global or even local optimality properties of biconvex optimization problems.
Instead, \emph{partial optimality} is usually considered for biconvex optimization problems, which is even weaker than local optimality.
Suppose $(x^\star, y^\star) \in \dset$ is a feasible point of (\ref{prob:bcp}), then $(x^\star, y^\star)$ is a \emph{partially optimal point} of (\ref{prob:bcp}) if
\begin{equation*}
    f_0(x^\star, y^\star) \leq f_0(x, y^\star)
\end{equation*}
for all $x \in \dset_{y^\star}$, and
\begin{equation*}
    f_0(x^\star, y^\star) \leq f_0(x^\star, y)
\end{equation*}
for all $y \in \dset_{x^\star}$, where
\begin{equation*}
    \dset_{y^\star} = \left\{x \in \xset\ \middle|\ \begin{array}{l}
        f_0(x, y^\star) < \infty\\
        f_i(x, y^\star) \leq 0,\quad i = 1, \ldots, m\\
        h_i(x, y^\star) = 0,\quad i = 1, \ldots, p
    \end{array}\right\},
\end{equation*}
and
\begin{equation*}
    \dset_{x^\star} = \left\{y \in \yset\ \middle|\ \begin{array}{l}
        f_0(x^\star, y) < \infty\\
        f_i(x^\star, y) \leq 0,\quad i = 1, \ldots, m\\
        h_i(x^\star, y) = 0,\quad i = 1, \ldots, p
    \end{array}\right\}.
\end{equation*}
It can be shown that for a differentiable biconvex optimization problem, every stationary point of $f_0$ over $\dset$ is partially optimal, and vice versa~\cite{gorski2007biconvex}.
However, such a point is not necessarily a local optimum, as stationary points can be saddle points of the objective function.
Some necessary conditions for a partially optimal point of a biconvex problem to be a local optimum are discussed in~\cite{gorski2007biconvex}, but in general, no stronger results can be obtained.
Although partial optimality is a weak notion of optimality, partially optimal points for biconvex problems are still widely used in practice, and often turn out to be good enough for many applications.

\section{Solving biconvex problems}\label{sec:bcp_solution}
As a result of the optimality properties regarding biconvex optimization problems discussed in \S\ref{sec:bcvxprob}, `solving' a biconvex optimization problem in practice usually amounts to finding a stationary point of the objective function over the feasible set.
In this section, we introduce the ACS heuristic for this purpose that is implemented in the DBCP framework, where the basic idea is to transform the biconvex problem into two convex subproblems, which can then be handled directly by DCP-type DSLs for convex optimization.
We also discuss some practical augmentations to ACS when dealing with the biconvex problem initialization and numerical stability issues.

\subsection{Alternate convex search}\label{sec:acs}
ACS as a minimization method is a special case of block relaxation methods~\cite{de1994block}, where the variables are divided into two disjoint blocks, and in each step, only the variables of an active block are optimized while those of the other block are fixed.
Specifically, for a biconvex optimization problem of the form (\ref{prob:bcp}), ACS iterates between solving the following two convex subproblems:
\begin{equation}\label{prob:bcp_x}
    \begin{array}{ll}
        \mbox{minimize} & f_0(x, \tilde{y})\\
        \mbox{subject to} & f_i(x, \tilde{y}) \leq 0,\quad i = 1, \ldots, m\\
        & h_i(x, \tilde{y}) = 0,\quad i = 1, \ldots, p,
    \end{array}
\end{equation}
where $x \in \xset$ is the variable, $\tilde{y} \in \yset$ is a fixed value of $y$, and
\begin{equation}\label{prob:bcp_y}
    \begin{array}{ll}
        \mbox{minimize} & f_0(\tilde{x}, y)\\
        \mbox{subject to} & f_i(\tilde{x}, y) \leq 0,\quad i = 1, \ldots, m\\
        & h_i(\tilde{x}, y) = 0,\quad i = 1, \ldots, p,
    \end{array}
\end{equation}
where $y \in \yset$ is the variable, and $\tilde{x} \in \xset$ is a fixed value of $x$.
Since the problems (\ref{prob:bcp_x}) and (\ref{prob:bcp_y}) are convex, efficient convex minimization methods can be used to solve these subproblems.
The full ACS procedure is summarized in the following algorithm.
\begin{algorithm}{Alternate convex search.}\label{alg:acs}
   \item \textbf{given} a starting point $(x^{(0)}, y^{(0)}) \in \dset$.
   \item $k \coloneqq 0$.
   \item \textbf{repeat}
      \begin{enumerate}
         \item \emph{Solve (\ref{prob:bcp_x}) with $\tilde{y} = y^{(k)}$.}
               Try to obtain
               \begin{equation}\label{eq:acs_x_update}
                  x^{(k+1)} \coloneqq \argmin_{x \in \xset} \left\{f_0(x, y^{(k)}) \biggm| \begin{array}{l}
                     f_i(x, y^{(k)}) \leq 0,\ i = 1, \ldots, m\\
                     h_i(x, y^{(k)}) = 0,\ i = 1, \ldots, p
                  \end{array}\right\},
               \end{equation}
               and quit if infeasible.
         \item \emph{Solve (\ref{prob:bcp_y}) with $\tilde{x} = x^{(k+1)}$.}
               Try to obtain
               \begin{equation}\label{eq:acs_y_update}
                  y^{(k+1)} \coloneqq \argmin_{y \in \yset} \left\{f_0(x^{(k+1)}, y) \biggm| \begin{array}{l}
                     f_i(x^{(k+1)}, y) \leq 0,\ i = 1, \ldots, m\\
                     h_i(x^{(k+1)}, y) = 0,\ i = 1, \ldots, p
                  \end{array}\right\},
               \end{equation}
               and quit if infeasible.
         \item $k \coloneqq k + 1$.
      \end{enumerate}
      \textbf{until} stopping criteria are satisfied.
\end{algorithm}

Gorski~\etal~\cite{gorski2007biconvex} showed that under weak assumptions, the sequence of objective function values ${\{f_0(x^{(k)}, y^{(k)})\}}_{k = 0}^\infty$ generated by ACS is monotonically nonincreasing and convergent.
Furthermore, if the sequence ${\{(x^{(k)}, y^{(k)})\}}_{k = 0}^\infty$ converges to $(x^\star, y^\star)$, then $(x^\star, y^\star)$ is a stationary point of $f_0$ (and hence is a partially optimal point of (\ref{prob:bcp})).
We should note that ACS is sensitive to initialization, \ie, with different initial points $(x^{(0)}, y^{(0)})$, the sequence ${\{(x^{(k)}, y^{(k)})\}}_{k = 0}^\infty$ generated through algorithm~\ref{alg:acs} may converge to different stationary points of $f_0$ with (probably) different function values.

There are several ways to define the stopping criteria for ACS.
A simple choice is to stop when the difference between the objective values of (\ref{prob:bcp}) obtained in two consecutive iterations is below a certain threshold $\epsilon > 0$, \ie, quit when
\begin{equation}\label{eq:stop_two_iter}
    |f_0(x^{(k+1)}, y^{(k+1)}) - f_0(x^{(k)}, y^{(k)})| < \epsilon.
\end{equation}
As a small variation, one may also use the difference between the optimal values of the problems (\ref{prob:bcp_x}) and (\ref{prob:bcp_y}) in one iteration as the criterion, \ie, quit when
\begin{equation}\label{eq:stop_one_iter}
    |f_0(x^{(k+1)}, y^{(k+1)}) - f_0(x^{(k+1)}, y^{(k)})| < \epsilon.
\end{equation}
In practice, there is not much difference between using (\ref{eq:stop_two_iter}) and (\ref{eq:stop_one_iter}) as the termination criterion for the ACS procedure, except that using the latter does not require storing the objective value of (\ref{prob:bcp}) from the previous iteration.
Note that when the objective values vary substantially in scale, it is useful in practice to combine the absolute threshold with a relative tolerance (see \S\ref{sec:impl} for this implementation).
Other choices include limiting the maximum number of iterations, or stopping when the changes in the optimization variables are below certain thresholds, \eg, quit when
\begin{equation*}
    \max\{\|x^{(k+1)} - x^{(k)}\|_2,\ \|y^{(k+1)} - y^{(k)}\|_2\} < \epsilon,
\end{equation*}
for some $\epsilon > 0$.

\subsection{Proximal regularizations}
Directly solving the subproblems (\ref{prob:bcp_x}) and (\ref{prob:bcp_y}) in the algorithm~\ref{alg:acs} may lead to numerical problems or slow convergence in practice.
One possible reason for this is that the biconvex objective function $f_0$ can be very `flat' in some region along certain directions when one block of variables is fixed (an example is shown in figure~\ref{fig:bcvx_fn}, in the region where $x$ and $y$ are close to zero), where the convex solver for subproblems may have difficulty making sufficient progress in minimizing the objective along the descent direction.
A common technique to alleviate these issues is to add proximal regularization terms to the objective functions of the subproblems.
Specifically, in the $(k+1)$th iteration of algorithm~\ref{alg:acs}, instead of performing (\ref{eq:acs_x_update}) and (\ref{eq:acs_y_update}), we try to perform the following updates:
\begin{equation}\label{eq:prox_acs_x_update}
    \begin{array}{llll}
        x^{(k+1)} &\coloneqq & \argmin_{x \in \xset} & f_0(x, y^{(k)}) + \lambda \norm{x - x^{(k)}}_2^2\\
        &&\mbox{subject to} & f_i(x, y^{(k)}) \leq 0,\quad i = 1, \ldots, m\\
        &&& h_i(x, y^{(k)}) = 0,\quad i = 1, \ldots, p,
    \end{array}
\end{equation}
and
\begin{equation}\label{eq:prox_acs_y_update}
    \begin{array}{llll}
        y^{(k+1)} &\coloneqq & \argmin_{y \in \yset} & f_0(x^{(k+1)}, y) + \lambda \norm{y - y^{(k)}}_2^2\\
        &&\mbox{subject to} & f_i(x^{(k+1)}, y) \leq 0,\quad i = 1, \ldots, m\\
        &&& h_i(x^{(k+1)}, y) = 0,\quad i = 1, \ldots, p,
    \end{array}
\end{equation}
where $\lambda \geq 0$ is the regularization parameter.
When $\lambda = 0$, the updates (\ref{eq:prox_acs_x_update}) and (\ref{eq:prox_acs_y_update}) reduce to (\ref{eq:acs_x_update}) and (\ref{eq:acs_y_update}), respectively.
When $\lambda > 0$, the proximal terms in (\ref{eq:prox_acs_x_update}) and (\ref{eq:prox_acs_y_update}) can be interpreted as adding soft trust region constraints to the respective optimization variables, which penalize large changes in the variables between two consecutive iterations, and in particular, the problems in (\ref{eq:prox_acs_x_update}) and (\ref{eq:prox_acs_y_update}) become strongly convex.
In practice, a large value of $\lambda$ helps improve numerical stability when solving the convex subproblems~\cite{boyd2011distributed,parikh2014proximal}.
It is also observed that adding the proximal regularizers can sometimes lead to better final points, \ie, points with lower objective values~\cite{shen2017disciplined}, compared to those from the original ACS procedure.
One possible intuition for this observation is that the proximal terms damp large changes in flat directions and favor nearby points when a subproblem has multiple or nearly equivalent solutions.
This changes the trajectory of the nonconvex ACS procedure and can lead to a different final point, although sensitivity to initialization remains.

\subsection{Initialization}
Note that algorithm~\ref{alg:acs} requires a feasible starting point $(x^{(0)}, y^{(0)}) \in \dset$, since otherwise, one of the subproblems (\ref{prob:bcp_x}) or (\ref{prob:bcp_y}) may be infeasible right at the first iteration.
Formally, the corresponding feasibility problem for a biconvex problem of the form (\ref{prob:bcp}) can be written as
\begin{equation}\label{prob:bcp_feas}
    \begin{array}{ll}
        \mbox{find} & (x, y)\\
        \mbox{subject to} & f_i(x, y) \leq 0,\quad i = 1, \ldots, m\\
        & h_i(x, y) = 0,\quad i = 1, \ldots, p
    \end{array}
\end{equation}
with variables $x \in \xset$ and $y \in \yset$.
Solving (\ref{prob:bcp_feas}) directly can be as hard as solving the original biconvex problem (\ref{prob:bcp}), and, actually, it can even be NP-hard~\cite{toker1995np}.
Here, we consider the following heuristic via relaxation to find a feasible starting point:
\begin{equation}\label{prob:bcp_feas_relax}
    \begin{array}{ll}
        \mbox{minimize} & \ones^T s + \norm{t}_1\\
        \mbox{subject to} & s \succeq 0\\
        & f_i(x, y) \leq s_i,\quad i = 1, \ldots, m\\
        & h_i(x, y) = t_i,\quad i = 1, \ldots, p,
    \end{array}
\end{equation}
where $x \in \xset$, $y \in \yset$ are the optimization variables, and $s \in \reals^m$ and $t \in \reals^p$ are the slack variables to relax the constraints in (\ref{prob:bcp_feas}).
The subscripts $i$ of $s_i$ and $t_i$ denote the $i$th entry of the vectors $s$ and $t$, respectively.
Note that the relaxed biconvex feasibility problem (\ref{prob:bcp_feas_relax}) is always feasible, since by choosing sufficiently large $s$ and $t$, all constraints can be satisfied.
To solve (\ref{prob:bcp_feas_relax}), we can again use the ACS procedure introduced in \S\ref{sec:acs}.
The full algorithm is given below.
\begin{algorithm}{Finding a feasible starting point.}\label{alg:bcp_init}
    \item \textbf{given} a starting point $(x^{(0)}, y^{(0)}) \in \xset \times \yset$.
    \item $k \coloneqq 0$.
    \item \textbf{repeat}
        \begin{enumerate}
            \item 
                $
                    \displaystyle{(x^{(k+1)}, s^\star, t^\star) \coloneqq \argmin_{\substack{x \in \xset\\s \in \reals^m, t \in \reals^p}} \left\{\ones^T s + \norm{t}_1\ \middle|\ \begin{array}{l}
                        s \succeq 0\\
                        f_i(x, y^{(k)}) \leq s_i,\ i = 1, \ldots, m\\
                        h_i(x, y^{(k)}) = t_i,\ i = 1, \ldots, p
                    \end{array}\right\}.}
                $

                \textbf{quit} with $(x^{(k+1)}, y^{(k)})$ \textbf{if} $\ones^T s^\star + \norm{t^\star}_1 = 0$.
            \item 
                $
                    \displaystyle{(y^{(k+1)}, s^\star, t^\star) \coloneqq \argmin_{\substack{y \in \yset\\s \in \reals^m, t \in \reals^p}} \left\{\ones^T s + \norm{t}_1\ \middle|\ \begin{array}{l}
                        s \succeq 0\\
                        f_i(x^{(k+1)}, y) \leq s_i,\ i = 1, \ldots, m\\
                        h_i(x^{(k+1)}, y) = t_i,\ i = 1, \ldots, p
                    \end{array}\right\}.}
                $

                \textbf{quit} with $(x^{(k+1)}, y^{(k+1)})$ \textbf{if} $\ones^T s^\star + \norm{t^\star}_1 = 0$.
            \item $k \coloneqq k + 1$.
        \end{enumerate}
        \textbf{until} maximum iterations are reached.
\end{algorithm}

Unlike algorithm~\ref{alg:acs}, algorithm~\ref{alg:bcp_init} requires only a point in $\xset \times \yset$, which is usually easy to obtain.
If algorithm~\ref{alg:bcp_init} quits with $\ones^T s^\star + \norm{t^\star}_1 = 0$ in some iteration, then the returned point is a feasible point of the original biconvex problem (\ref{prob:bcp}), which can then be used as a starting point for the ACS procedure in algorithm~\ref{alg:acs}.
However, we must note that there is no guarantee that the algorithm~\ref{alg:bcp_init} will find a feasible point for any instance of the biconvex problem (\ref{prob:bcp}), even if such a point exists.
In practice, as a generic practical method, this approach seems to work quite well.

\subsection{Penalty method and infeasible start}
Now we integrate the relaxation, which transforms the biconvex feasibility problem (\ref{prob:bcp_feas}) into (\ref{prob:bcp_feas_relax}), directly into the original biconvex problem (\ref{prob:bcp}), so that the ACS procedure can be applied even when starting from an infeasible point.
Let $s \in \reals^m$ and $t \in \reals^p$ be slack variables to relax the inequality and equality constraints of (\ref{prob:bcp}), respectively.
We consider the following relaxed biconvex optimization problem:
\begin{equation}\label{prob:bcp_relax}
    \begin{array}{ll}
        \mbox{minimize} & f_0(x, y) + \nu(\ones^T s + \norm{t}_1)\\
        \mbox{subject to} & s \succeq 0\\
        & f_i(x, y) \leq s_i,\quad i = 1, \ldots, m\\
        & h_i(x, y) = t_i,\quad i = 1, \ldots, p,
    \end{array}
\end{equation}
where $\nu > 0$ is a penalty parameter.
The problem (\ref{prob:bcp_relax}) is always feasible, since by choosing sufficiently large $s$ and $t$, all constraints can be satisfied.
Moreover, if the original biconvex problem (\ref{prob:bcp}) is feasible, then for sufficiently large $\nu$, applying ACS to (\ref{prob:bcp_relax}) will yield a final point $(x^\star, y^\star, s^\star, t^\star)$, such that $\ones^T s^\star + \norm{t^\star}_1 = 0$, \ie, $(x^\star, y^\star)$ is feasible and partially optimal for (\ref{prob:bcp})~\cite{nocedal2006numerical,shen2017disciplined}.
The full algorithm is given as follows.
\begin{algorithm}{Penalty alternate convex search.}\label{alg:infeas_acs}
    \item \textbf{given} a starting point $(x^{(0)}, y^{(0)}) \in \xset \times \yset$ and sufficiently large $\nu > 0$.
    \item $k \coloneqq 0$.
    \item \textbf{repeat}
        \begin{enumerate}
            \item 
                $
                    \displaystyle{x^{(k+1)} \coloneqq \argmin_{x \in \xset} \left\{\begin{array}{l}
                        f_0(x, y^{(k)})\\
                        \quad+ \nu(\ones^T s + \norm{t}_1)
                    \end{array}\ \middle|\ \begin{array}{l}
                        s \in \reals^m,\quad t \in \reals^p,\quad s \succeq 0\\
                        f_i(x, y^{(k)}) \leq s_i,\quad i = 1, \ldots, m\\
                        h_i(x, y^{(k)}) = t_i,\quad i = 1, \ldots, p
                    \end{array}\right\}.}
                $
            \item 
                $
                    \displaystyle{y^{(k+1)} \coloneqq \argmin_{y \in \yset} \left\{\begin{array}{l}
                        f_0(x^{(k+1)}, y)\\
                        \quad+ \nu(\ones^T s + \norm{t}_1)
                    \end{array}\ \middle|\ \begin{array}{l}
                        s \in \reals^m,\quad t \in \reals^p,\quad s \succeq 0\\
                        f_i(x^{(k+1)}, y) \leq s_i,\quad i = 1, \ldots, m\\
                        h_i(x^{(k+1)}, y) = t_i,\quad i = 1, \ldots, p
                    \end{array}\right\}.}
                $
            \item $k \coloneqq k + 1$.
        \end{enumerate}
        \textbf{until} stopping criteria are satisfied.
\end{algorithm}

Compared to algorithm~\ref{alg:acs}, the penalty ACS procedure in algorithm~\ref{alg:infeas_acs} can start from any point in $\xset \times \yset$, and hence, the initialization step via, \eg, algorithm~\ref{alg:bcp_init}, can be avoided.
The same termination criteria as those discussed in \S\ref{sec:acs} can still be used for algorithm~\ref{alg:infeas_acs}.
However, we must note that there is no guarantee that the final point returned by algorithm~\ref{alg:infeas_acs} is feasible for the original biconvex problem (\ref{prob:bcp}), even if such a point exists, since the penalty parameter $\nu$ may not be sufficiently large.
In practice, the value of $\nu$ can be selected in an ad hoc manner, \ie, one may try to increase $\nu$ and resolve (\ref{prob:bcp_relax}) if the final point returned by algorithm~\ref{alg:infeas_acs} is still infeasible for (\ref{prob:bcp}).
Finally, the proximal regularizations as in (\ref{eq:prox_acs_x_update}) and (\ref{eq:prox_acs_y_update}) can be readily integrated into the subproblems in the algorithm~\ref{alg:infeas_acs} to improve numerical stability.

\section{Disciplined biconvex programming}\label{sec:dbcp}
We now present the DBCP biconvex ruleset for modeling biconvex optimization problems in a way that the biconvexity is easily verified by construction.

\paragraph{DCP ruleset.}
The DBCP ruleset is heavily based on the DCP ruleset, which consists of a library of convex atomic functions and a convex syntax ruleset that prescribes how these atomic functions may be composed to express (more complex) convex optimization problems~\cite{grant2006disciplined}.
Specifically, all functions in a DCP-compliant problem must be formed as an \emph{expression} consisting of variables, constants or parameters, and atomic functions.
The sign and curvature of each DCP expression can be determined recursively using the signs and curvatures of its arguments together with the sign, curvature, and possibly sign-dependent monotonicity of the atomic function applied at each step. Based on this information, the convexity of the overall problem can be verified by checking whether each function in the problem satisfies the DCP composition rules~\cite{grant2006disciplined,boyd2004convex}.
We make an observation that is critical for extending DCP to DBCP:
The basic arithmetic properties, in particular, the composition property, of biconvex functions, as discussed in \S\ref{sec:bcvx_fn}, are compatible with the DCP ruleset.
In other words, an optimization problem formed by verifiable biconvex expressions is still biconvex if it complies with a variant of the DCP ruleset in which `convex' is replaced by `biconvex'.

\paragraph{DBCP product rule.}
To construct disciplined biconvex expressions, we first note that the product of expressions that are both nonconstant is prohibited in DCP, since the convexity of such expressions cannot be determined in general~\cite{grant2006disciplined};
however, most biconvex expressions that appear in practice are constructed through variable multiplications.
Hence, we introduce the following product rule for DBCP:
\begin{itemize}
    \item
        Products involving a constant or parameter factor are valid in DBCP whenever they satisfy the DCP rules. If both the left-hand and right-hand expressions contain variables, a valid DBCP biconvex product expression should be one of the following forms:
        \begin{center}
            \emph{affine} \texttt{*} \emph{affine}\\
            \emph{affine-nonneg} \texttt{*} \emph{convex}\quad or\quad \emph{affine-nonpos} \texttt{*} \emph{concave}\\
            \emph{convex-nonneg} \texttt{*} \emph{convex-nonneg}\quad or\quad \emph{concave-nonpos} \texttt{*} \emph{concave-nonpos}
        \end{center}
        The \emph{nonneg} and \emph{nonpos} qualifiers indicate that the expression is known to be nonnegative or nonpositive, respectively.
    \item
        Across the objective and constraints, the same two-block variable assignment must be used for every product in which both factors contain variables. For each such product, one factor must contain only variables from one block and the other only variables from the other block.
\end{itemize}

The common two-block assignment requirement in the last bullet prevents expressions like \texttt{x * y}, \texttt{y * z}, and \texttt{z * x} from appearing simultaneously in the same optimization problem, since \texttt{x * y} and \texttt{y * z} require \texttt{x} and \texttt{z} to belong to the same block, whereas \texttt{z * x} requires them to belong to different blocks.
According to the definition and basic properties of biconvex functions listed in \S\ref{sec:bcvx_fn}, expressions formed according to the above DBCP product rule are guaranteed to be biconvex.

\section{Implementation}\label{sec:impl}
In this section, we introduce our Python implementation of the DBCP framework, as an extension to CVXPY~\cite{diamond2016cvxpy,agrawal2018rewriting}.
The corresponding open-source package, named \texttt{dbcp}, is available at
\begin{quote}
    \url{https://github.com/dxogrp/dbcp}.
\end{quote}

\paragraph{Biconvex atoms.}
Biconvex expressions in \texttt{dbcp} are created from fundamental atomic functions.
Most of the atoms used in \texttt{dbcp} are inherited from CVXPY, such as inner product and elementwise multiplication.
An example of the extended atoms that are not supported by CVXPY is the one-dimensional discrete \emph{convolution} operation between two one-dimensional nonconstant expressions, which can be called by the user via \texttt{dbcp.convolve}.
It has the same behavior as the \texttt{convolve} function in CVXPY, except that the first expression need not be a constant.
The biconvex atom library is extensible, so that additional atoms can be added as necessary.

\paragraph{Specifying biconvex problems.}
Users can define their optimization variables, objective functions, and constraints using the standard CVXPY syntax, and a biconvex problem is constructed using:
\begin{quote}
    \texttt{prob = BiconvexProblem(obj, x\_var, y\_var, constr)}
\end{quote}
The argument \texttt{obj} is the objective constructed from a DBCP-compliant biconvex expression, \texttt{x\_var} and \texttt{y\_var} are iterables of the biconvex optimization variables, and \texttt{constr} is a list of DBCP-compliant biconvex constraints.
The arguments \texttt{x\_var} and \texttt{y\_var} define the two disjoint groups of variables that are fixed on alternating steps during the ACS procedure, so that the two subproblems are convex.
Note that \texttt{x\_var} and \texttt{y\_var} usually contain all the problem variables as in (\ref{prob:bcp}), but this is not necessarily the case; our implementation allows for more flexibility in defining these groups.
Specifically, problem variables omitted from both \texttt{x\_var} and \texttt{y\_var} can remain active in both subproblems, provided both resulting subproblems are DCP-compliant.
For example, $Z$ is such a variable in the problem
\begin{equation*}
    \begin{array}{ll}
        \mbox{minimize} & \norm{XY + Z - A}_F\\
        \mbox{subject to} & \norm{Z}_F \leq 1,
    \end{array}
\end{equation*}
where $X \in \reals^{m \times k}$, $Y \in \reals^{k \times n}$, $Z \in \reals^{m \times n}$ are optimization variables, and $A \in \reals^{m \times n}$ is problem data.
To specify this biconvex problem, one may use the following code:
\begin{lstlisting}[language=zhpython]
import cvxpy as cp
from dbcp import BiconvexProblem

X = cp.Variable((m, k))
Y = cp.Variable((k, n))
Z = cp.Variable((m, n))

obj = cp.Minimize(cp.norm(X @ Y + Z - A, 'fro'))
constr = [cp.norm(Z, 'fro') <= 1]
prob = BiconvexProblem(obj, [X], [Y], constr)
\end{lstlisting}

\paragraph{Solving a biconvex problem.}
To solve the specified biconvex problem \texttt{prob}, one may simply call \texttt{prob.solve()}.
There are many optional arguments that can be passed to \texttt{prob.solve()} to customize the solution procedure, most of which are directly inherited from CVXPY for the configuration of the convex solvers used to solve the subproblems.
In addition to these standard arguments, \texttt{dbcp} also provides several specific options for controlling the ACS procedure.
One of the most important ones is \texttt{lbd}, which sets a value for the proximal regularization parameter $\lambda$ in (\ref{eq:prox_acs_x_update}) and (\ref{eq:prox_acs_y_update}).
When \texttt{lbd} is zero, no proximal regularization is added, \ie, the original ACS procedure (algorithm~\ref{alg:acs}) is used;
if \texttt{lbd} is positive, then the proximal regularized updates in (\ref{eq:prox_acs_x_update}) and (\ref{eq:prox_acs_y_update}) are used instead of the original updates (\ref{eq:acs_x_update}) and (\ref{eq:acs_y_update}).
It is recommended that the user specify (feasible) initial values for all optimization variables before calling \texttt{prob.solve()}.
Otherwise, \texttt{dbcp} will try to generate random initial values for the unspecified variables from a standard normal distribution, which is not guaranteed to work well for all problems.
With the default setting \texttt{mode='direct'}, the \texttt{prob.solve()} method first checks the feasibility of the current initial point.
If the initial point is infeasible, then algorithm~\ref{alg:bcp_init} will be used to try to find a feasible starting point, based on the current values of the variables, before launching the ACS procedure.
We implement a scale-aware version of the stopping criterion (\ref{eq:stop_one_iter}).
Let $u = f_0(x^{(k+1)}, y^{(k)})$ and $v = f_0(x^{(k+1)}, y^{(k+1)})$ denote the objective values of the two subproblems in one iteration.
The ACS procedure terminates when
\begin{equation}\label{eq:scale_aware_stop}
    |u-v| \leq \epsilon_{\mathrm{abs}} + \epsilon_{\mathrm{rel}}\max\{|u|,|v|\}.
\end{equation}
The argument \texttt{abs\_tol} to \texttt{prob.solve()} sets $\epsilon_{\mathrm{abs}}$ and defaults to $10^{-6}$.
The argument \texttt{rel\_tol} sets $\epsilon_{\mathrm{rel}}$ and also defaults to $10^{-6}$.
Setting \texttt{rel\_tol} to zero gives a criterion involving only the absolute tolerance.

\paragraph{Selecting the penalty method.}
Setting \texttt{mode='penalty'} in \texttt{prob.solve()} directly solves the relaxed biconvex problem (\ref{prob:bcp_relax}) using the penalty ACS procedure (algorithm~\ref{alg:infeas_acs}) from the current variable values, thereby permitting an infeasible starting point.
If initial values are not specified, random values are generated as with the direct method.
The option \texttt{nu} sets the positive penalty parameter $\nu$ in (\ref{prob:bcp_relax}), and \texttt{slack\_tol}, which defaults to $10^{-6}$, sets the tolerance for the final total slack.
Note that although the problem (\ref{prob:bcp_relax}) is always feasible, there is no guarantee that the final point returned by \texttt{prob.solve()} is feasible for the original biconvex problem (\ref{prob:bcp}).
To monitor the feasibility progress of the algorithm, the value of the total slack $\ones^T s + \norm{t}_1$ is reported after each ACS iteration.
In the case that the final point is still infeasible with total slack exceeding \texttt{slack\_tol}, the user may try to increase the value of \texttt{nu} and resolve the same problem object.

\paragraph{Solution statuses.}
After the ACS procedure terminates normally, the DBCP termination status is stored in \texttt{prob.status} and takes one of the following four values:
\begin{itemize}
    \item \texttt{converged}:
        The stopping criterion (\ref{eq:scale_aware_stop}) is satisfied, and if the penalty method is used, the final total slack does not exceed \texttt{slack\_tol}.
    \item \texttt{converged\_with\_slack}:
        The stopping criterion (\ref{eq:scale_aware_stop}) is satisfied, but the final total slack exceeds \texttt{slack\_tol}.
    \item \texttt{converged\_inaccurate}:
        The specified maximum number of ACS iterations is reached before the stopping criterion is satisfied, and if the penalty method is used, the final total slack does not exceed \texttt{slack\_tol}.
    \item \texttt{converged\_inaccurate\_with\_slack}:
        The maximum number of ACS iterations is reached before the stopping criterion is satisfied, and the final total slack exceeds \texttt{slack\_tol}.
\end{itemize}
Note that for the direct method, \texttt{prob.status} can only take the two values without the \texttt{\_with\_slack} suffix.

\paragraph{Verification of biconvexity.}
The \texttt{dbcp} package automatically verifies the biconvexity of the objective function and all constraints according to the DBCP ruleset introduced in \S\ref{sec:dbcp} when a biconvex optimization problem is specified as a \texttt{BiconvexProblem} instance.
The user can check whether a problem is DBCP-compliant by calling \texttt{prob.is\_dbcp()}, which returns \texttt{True} if the problem is a valid DBCP biconvex problem, and \texttt{False} otherwise.
If a user is trying to solve a non-DBCP-compliant biconvex problem, \texttt{dbcp} will raise an error directly after the \texttt{prob.solve()} method is called.

\paragraph{Generalized inequality constraints.}
Lastly, we should mention that although all the previous discussion about solving and finding a feasible initial point for biconvex problems is based on standard inequality constraints as shown in (\ref{prob:bcp}), the \texttt{dbcp} package is implemented in a way that natively supports biconvex problems with generalized inequality constraints, such as second-order cone constraints and positive semidefinite constraints.
Interested readers may refer to \S\ref{sec:gen_ineq_constr} for more theoretical details about how these constraints are handled in the backend, while practitioners can directly use the standard CVXPY syntax for specifying generalized inequality constraints when defining biconvex problems using \texttt{dbcp}.
The appearance of such constraints will be automatically detected by \texttt{dbcp}, and the appropriate adaptations will be made to the solution procedure without any extra effort from the user.

\section{Examples}\label{sec:example}
In this section, we present several numerical examples using the \texttt{dbcp} package to specify and solve biconvex optimization problems that appear frequently in practice.
We only show snippets of the code for defining and solving the biconvex problems, while omitting the data generation and result visualization parts for brevity.
Interested readers can find the full code examples at
\begin{quote}
    \url{https://github.com/dxogrp/dbcp/tree/main/examples}.
\end{quote}

\subsection{Nonnegative matrix factorization}\label{sec:eg_nmf}
We start with a basic nonnegative matrix factorization problem.
Suppose that we are given a matrix $A \in \reals^{m \times n}$, and are interested in finding two nonnegative matrices $X \in \reals^{m \times k}$ and $Y \in \reals^{k \times n}$ such that $A \approx XY$.
This can be formulated as the following biconvex optimization problem:
\begin{equation*}
    \begin{array}{ll}
        \mbox{minimize} & \norm{XY - A}_F^2\\
        \mbox{subject to} & X_{ij} \geq 0,\quad i = 1, \ldots, m,\quad j = 1, \ldots, k\\
        & Y_{ij} \geq 0,\quad i = 1, \ldots, k,\quad j = 1, \ldots, n
    \end{array}
\end{equation*}
with variables $X$ and $Y$.
To specify this problem using \texttt{dbcp}, one may use the following code:
\begin{lstlisting}[language=zhpython]
X = cp.Variable((m, k), nonneg=True)
Y = cp.Variable((k, n), nonneg=True)

obj = cp.Minimize(cp.sum_squares(X @ Y - A))
prob = BiconvexProblem(obj, [X], [Y])
\end{lstlisting}

In this example, we set $m = 5$, $n = 10$, and $k = 5$, generate two factor matrices with entries sampled uniformly from $[0,1)$, and form $A$ as their product.
After calling \texttt{prob.solve()}, we obtain a final point with an objective value around $6.09 \times 10^{-6}$, which indicates that the original matrix $A$ is well approximated by the product of the recovered nonnegative matrices $X$ and $Y$.

\subsection{Bilinear logistic regression}
We consider a bilinear logistic regression problem for binary classification~\cite{dyrholm2007bilinear}.
Suppose that we are given a dataset $(X_i, y_i)$, $i = 1, \ldots, m$, where each sample consists of a feature matrix $X_i \in \reals^{n \times k}$ and a binary label $y_i \in \{0, 1\}$.
Our goal is to construct a bilinear classifier $\hat{y} = 1$ if $\tr(U^T X V) > 0$, and $\hat{y} = 0$ otherwise, where $U \in \reals^{n \times r}$ and $V \in \reals^{k \times r}$ are the bilinear logistic regression coefficients with a predefined (maximum) rank $r$, and $\tr(M)$ denotes the trace of some square matrix $M$.
To fit a bilinear logistic regression model to the dataset, we would like to solve the following bilinear maximum likelihood estimation problem:
\begin{equation*}
    \begin{array}{ll}
        \mbox{maximize} & \sum_{i = 1}^{m} y_i \tr(U^T X_i V) - \log(1 + \exp(\tr(U^T X_i V)))
    \end{array}
\end{equation*}
with variables $U$ and $V$.
This problem can be specified using \texttt{dbcp} as follows:
\begin{lstlisting}[language=zhpython]
U = cp.Variable((n, r))
V = cp.Variable((k, r))

obj = cp.Maximize(sum(
    y * cp.trace(U.T @ X @ V)
    - cp.logistic(cp.trace(U.T @ X @ V))
    for X, y in zip(Xs, ys)))
prob = BiconvexProblem(obj, [U], [V])
\end{lstlisting}

In this example, we set $m = 300$, $n = 20$, $k = 10$, and $r = 5$.
The dataset $(X_i, y_i)$, $i = 1, \ldots, m$, is generated synthetically via the \texttt{make\_classification} function from the scikit-learn library~\cite{pedregosa2011scikit}.
After calling \texttt{prob.solve()} with the proximal regularization weight \texttt{lbd=0.1} and the termination threshold \texttt{abs\_tol=1e-4}, we obtain a final point with an objective value $-1.44 \times 10^{-3}$.

\subsection{$k$-means clustering}
Suppose that we are given a set of data points $x_i \in \reals^n$, $i = 1, \ldots, m$, and we would like to cluster them into $k$ groups, using the $k$-means clustering method.
This problem can be formulated as the following biconvex optimization problem~\cite{zhu2025multi}:
\begin{equation}\label{prob:kmeans}
    \begin{array}{ll}
        \mbox{minimize} & \sum_{i = 1}^{m} z_i^T (\norm{\bar{x}_1 - x_i}_2^2, \ldots, \norm{\bar{x}_k - x_i}_2^2)\\
        \mbox{subject to} & 0 \preceq z_i \preceq \ones,\quad \ones^T z_i = 1,\quad i = 1, \ldots, m
    \end{array}
\end{equation}
with variables $\bar{x}_i \in \reals^n$, $i = 1, \ldots, k$, and $z_i \in \reals^k$, $i = 1, \ldots, m$.

We can interpret the variables in the problem (\ref{prob:kmeans}) as follows:
The variables $\bar{x}_1, \ldots, \bar{x}_k$ represent the cluster centroids, and each variable $z_i$ is a soft assignment vector for the data point $x_i$, where the $j$th entry of $z_i$ indicates the probability that the sample $x_i$ belongs to cluster $j$.
Then, the objective function in (\ref{prob:kmeans}) represents the total within-cluster variance, which we would like to minimize.
This formulation relaxes the original $k$-means problem by allowing each $z_i$ to be a soft rather than a hard assignment.
However, the relaxation is tight, since for fixed centroids, the objective is minimized by assigning each data point entirely to a nearest centroid (provided there is no tie in the distances).
Therefore, the relaxed problem (\ref{prob:kmeans}) always admits a hard assignment optimum and has the same optimal value as the original $k$-means problem.

To specify the problem (\ref{prob:kmeans}) using \texttt{dbcp}, one may use the following code:
\begin{lstlisting}[language=zhpython]
xbars = cp.Variable((k, n))
zs = cp.Variable((m, k), nonneg=True)

d = cp.sum_squares(xs[:, None, :] - xbars[None, :, :], axis=2)
obj = cp.Minimize(cp.sum(cp.multiply(zs, d)))
constr = [zs <= 1, cp.sum(zs, axis=1) == 1]
prob = BiconvexProblem(obj, [xbars], [zs], constr)
\end{lstlisting}

\begin{figure}[t]
    \centering
    \includegraphics[width=0.6\textwidth]{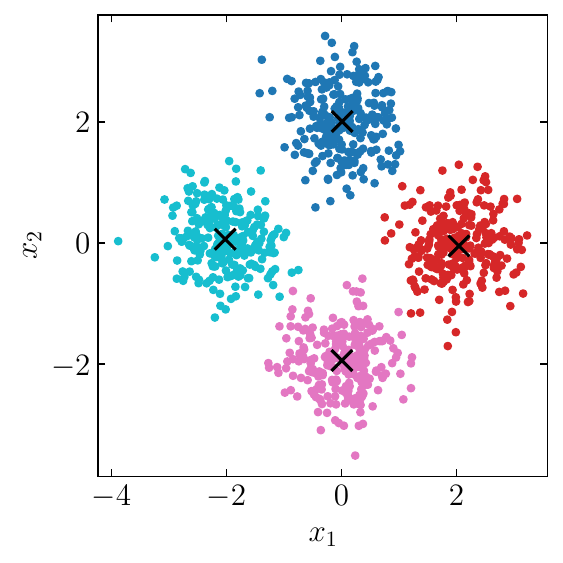}%
    \caption{Results of the $k$-means clustering example.}\label{fig:kmeans}
\end{figure}

We generate a synthetic dataset of $m = 1000$ points in $\reals^2$ using the \texttt{make\_blobs} function of scikit-learn, and set the number of clusters $k = 4$ with the ground truth centroids $(0, 2)$, $(0, -2)$, $(2, 0)$, and $(-2, 0)$.
The $k$-means clustering results based on the formulation (\ref{prob:kmeans}) are shown in figure~\ref{fig:kmeans}.
After \texttt{prob.solve()} terminates, the displayed hard cluster label $c_i$ of each data point is obtained from the optimized soft assignment vector $z_i$ according to
\begin{equation*}
    c_i = \argmax_{j \in \{1, \ldots, k\}} z_{ij},\quad i = 1, \ldots, m,
\end{equation*}
where the ties are resolved by choosing the smallest maximizing index.
The colors indicate these cluster labels $c_1, \ldots, c_m$, and the black crosses represent the optimized centroids $\bar{x}_1, \ldots, \bar{x}_k$.

\subsection{Dictionary learning}
We consider the sparse dictionary learning problem~\cite{aharon2006k}, which aims to find a dictionary matrix $D \in \reals^{m \times k}$ and a sparse code matrix $X \in \reals^{k \times n}$, such that the data matrix $Y \in \reals^{m \times n}$ can be well approximated by their product $DX$, while the matrix $X$ is sparse and the matrix $D$ has bounded Frobenius norm.
The dictionary learning problem can be formulated as the following biconvex optimization problem:
\begin{equation}\label{prob:dict_learning}
    \begin{array}{ll}
        \mbox{minimize} & \norm{DX - Y}_F^2 + \alpha \norm{X}_1\\
        \mbox{subject to} & \norm{D}_F \leq \beta
    \end{array}
\end{equation}
with variables $D$ and $X$, where $\alpha > 0$ is the sparsity regularization parameter, and $\beta > 0$ is the bound on the Frobenius norm of the dictionary matrix.
To specify this problem using \texttt{dbcp}, one may use the following code:
\begin{lstlisting}[language=zhpython]
D = cp.Variable((m, k))
X = cp.Variable((k, n))

obj = cp.Minimize(
    cp.sum_squares(D @ X - Y) + alpha * cp.norm1(X)
)
prob = BiconvexProblem(
    obj, [D], [X], [cp.norm(D, 'fro') <= beta]
)
\end{lstlisting}

\begin{figure}[t]
    \centering
    \includegraphics[width=0.8\textwidth]{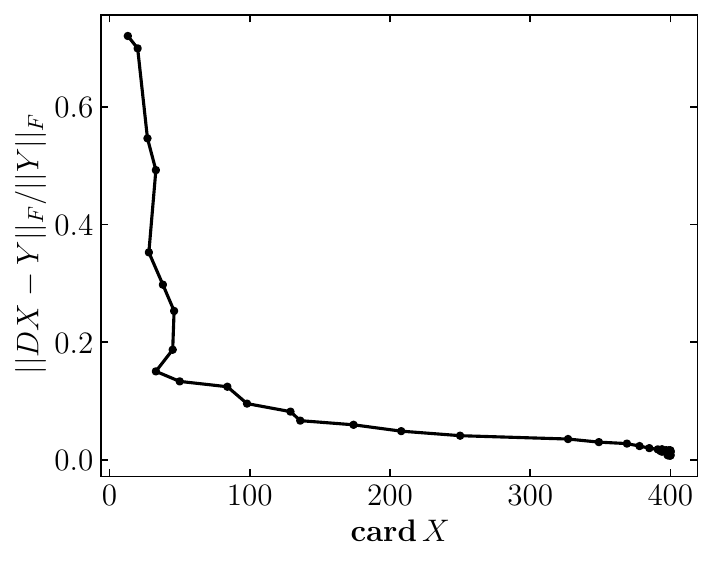}%
    \caption{Trade-off curve of the sparse dictionary learning example.}\label{fig:dict_learning}
\end{figure}

In this example, we set $m = 10$, $n = 20$, $k = 20$, $\beta = 1$, and generate the data matrix $Y \in \reals^{m \times n}$ from the standard normal distribution.
The problem (\ref{prob:dict_learning}) is then solved for different values of the sparsity regularization parameter $\alpha$ ranging from $10^{-5}$ to $1$.
For each value of $\alpha$, the values of the problem variables $D$ and $X$ are reinitialized randomly before \texttt{prob.solve()} is called.
Figure~\ref{fig:dict_learning} shows the trade-off curve between the relative approximation error $\norm{DX-Y}_F / \norm{Y}_F$ and the cardinality (\ie, the number of nonzero entries) of $X$.
We should note that this curve is not a Pareto frontier, as some plotted points are Pareto-dominated by others.
The ACS procedure is a local method and need not return a global solution for each value of $\alpha$.
Moreover, the returned point can depend on the random initialization.
The optimization problem uses $\norm{X}_1$ as a surrogate for sparsity, whereas the horizontal axis shows the cardinality of $X$, computed here as the number of entries whose absolute values are at least $10^{-3}$.
Since the $\ell_1$-norm also depends on the magnitudes of the nonzero entries, even a global solution of problem~(\ref{prob:dict_learning}) need not be Pareto optimal in terms of approximation error and cardinality.

\subsection{Blind deconvolution}
Blind deconvolution is a technique used to recover some sharp signal or image from a blurred observation when the blur itself is unknown~\cite{chaudhuri2014blind}.
It jointly estimates the original signal and the blur kernel, with some prior knowledge about their structures.
Suppose that we are given a data vector $d \in \reals^{m + n - 1}$, which is the convolution of an unknown sparse signal $x \in \reals^n$ and an unknown smooth vector $y \in \reals^m$ with bounded $\ell_\infty$-norm (\ie, bounded largest absolute entry).
Additionally, we have the prior knowledge that both vectors $x$ and $y$ are nonnegative.
The corresponding blind deconvolution problem can be formulated as the following biconvex optimization problem:
\begin{equation*}
    \begin{array}{ll}
        \mbox{minimize} & \norm{x \otimes y - d}_2^2 + \alpha_{\rm sp} \norm{x}_1 + \alpha_{\rm sm} \norm{Dy}_2^2\\
        \mbox{subject to} & x \succeq 0,\quad y \succeq 0\\
        & \norm{y}_\infty \leq \beta
    \end{array}
\end{equation*}
with variables $x$ and $y$, where $\alpha_{\rm sp}, \alpha_{\rm sm} > 0$ are the regularization parameters for the sparsity of $x$ and smoothness of $y$, respectively, and $\beta > 0$ is the bound on the $\ell_\infty$-norm of the vector $y$.
The matrix $D \in \reals^{(m - 1) \times m}$ is the first-order difference operator, given by
\begin{equation*}
    D = \left[\begin{array}{ccccc}
        1 & -1 &&&\\
        & 1 & -1 &&\\
        && \ddots & \ddots &\\
        &&& 1 & -1
    \end{array}\right],
\end{equation*}
so that $Dy$ computes the vector of successive differences of $y$.
The convolution $x \otimes y$ of the vectors $x$ and $y$ is given by
\begin{equation*}
    {(x \otimes y)}_k = \sum_{i + j = k + 1} x_i y_j,\quad k = 1, \ldots, m + n - 1.
\end{equation*}
To specify this problem using \texttt{dbcp}, one may use the following code:
\begin{lstlisting}[language=zhpython]
x = cp.Variable(n, nonneg=True)
y = cp.Variable(m, nonneg=True)

obj = cp.Minimize(
    cp.sum_squares(convolve(x, y) - d)
    + alpha_sp * cp.norm1(x)
    + alpha_sm * cp.sum_squares(cp.diff(y)))
constr = [cp.norm(y, 'inf') <= beta]
prob = BiconvexProblem(obj, [x], [y], constr)
\end{lstlisting}

\begin{figure}[t]
    \centering
    \includegraphics[width=0.8\textwidth]{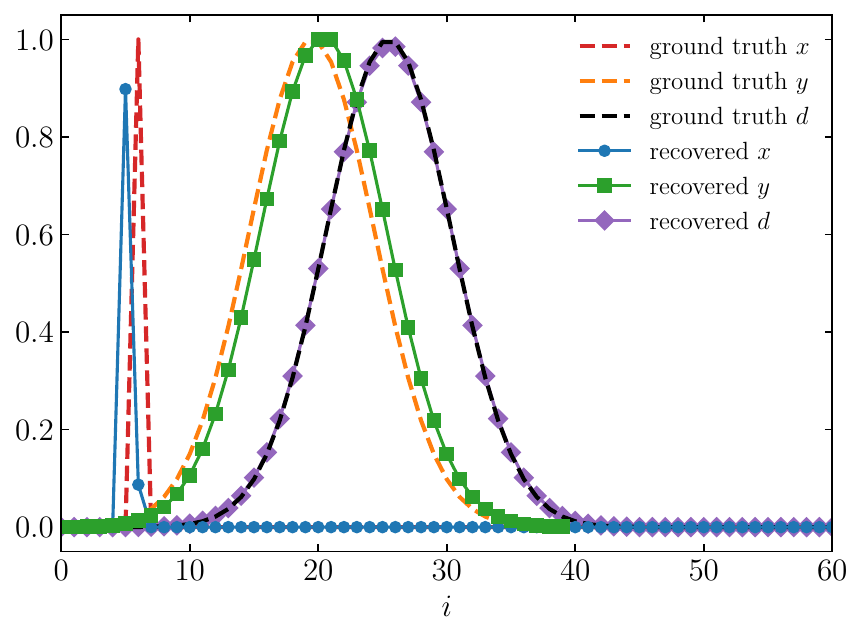}%
    \caption{Results of the blind deconvolution example.}\label{fig:blind_deconv}
\end{figure}

In this example, we consider the vectors $x \in \reals^{120}$, $y \in \reals^{40}$, and the $\ell_\infty$-norm bound $\beta = 1$.
The ground truth vectors were generated according to a similar example used by Shen~\etal~\cite{shen2017disciplined}.
The ground truth and the recovered signals after calling \texttt{prob.solve()} with regularization parameters $\alpha_{\rm sp} = 0.1$ and $\alpha_{\rm sm} = 0.2$ are shown in figure~\ref{fig:blind_deconv}.

\subsection{Fitting input-output hidden Markov models}
We consider the problem of fitting a logistic input-output hidden Markov model (IO-HMM) to a dataset~\cite{jha2024active}.
Suppose that we are given a dataset $(x(t), y(t))$, $t = 1, \ldots, m$, where each sample consists of an input feature vector $x(t) \in \reals^n$ and an output label $y(t) \in \{0, 1\}$, generated from a $K$-state IO-HMM, according to the following procedure:
Let $\hat{z}(t) \in \{1, \ldots, K\}$, $t = 1, \ldots, m$, be the state label of the IO-HMM with initial state distribution $p_{\rm init} \in \reals^K$ with $\ones^T p_{\rm init} = 1$ and transition matrix $P_{\rm tr} \in \reals^{K \times K}$ with $P_{\rm tr} \ones = \ones$.
At the time step $t$, the state label $\hat{z}(t)$ is sampled according to
\begin{equation*}
    \hat{z}(t) \sim \left\{
        \begin{array}{ll}
            {\rm Cat}(p_{\rm init}), & t = 1\\
            {\rm Cat}(p_{\hat{z}(t - 1)}), & t > 1,
        \end{array}\right.
\end{equation*}
where the vector $p_{\hat{z}(t-1)} \in \reals^K$ denotes the $\hat{z}(t-1)$th row of the matrix $P_{\rm tr}$, and ${\rm Cat}(p)$ denotes the categorical distribution with $p$ being the vector of category probabilities.
Then, given the feature vector $x(t) \in \reals^n$, the output $y(t) \in \{0, 1\}$ of this IO-HMM at time step $t$ is generated from a logistic model, \ie,
\begin{equation*}
    \prob(y(t) = 1) = \frac{1}{1 + \exp(-{x(t)}^T \theta_{\hat{z}(t)})},
\end{equation*}
where $\theta_{\hat{z}(t)} \in \{\theta_1, \ldots, \theta_K\} \subseteq \reals^n$ is the coefficient.

Given the dataset $(x(t), y(t))$, $t = 1, \ldots, m$, we are interested in recovering the transition matrix $P_{\rm tr}$, the model parameters $\theta_1, \ldots, \theta_K$, and the unobserved state labels $\hat{z}(1), \ldots, \hat{z}(m)$.
For each $t = 1, \ldots, m$, define the vector-valued expression $r(t) = (r_1(t), \ldots, r_K(t)) \in \reals^K$ componentwise by
\begin{equation*}
    r_k(t) = y(t){x(t)}^T \theta_k - \log(1 + \exp({x(t)}^T \theta_k)),\quad k = 1, \ldots, K.
\end{equation*}
Noticing that the transition matrix $P_{\rm tr}$ can be easily estimated from the state labels $\hat{z}(t)$, $t = 1, \ldots, m$, we consider the following biconvex optimization problem for fitting the IO-HMM~\cite{zhu2025multi}:
\begin{equation}\label{prob:iohmm}
    \begin{array}{ll}
        \mbox{minimize} & -\sum_{t = 1}^{m} {z(t)}^T r(t)\\
        &\qquad + \alpha_\theta \sum_{k = 1}^{K} \norm{\theta_k}^2_2 + \alpha_z \sum_{t = 1}^{m - 1} \Dkl(z(t), z(t + 1))\\
        \mbox{subject to} & 0 \preceq z(t) \preceq \ones,\quad \ones^T z(t) = 1,\quad t = 1, \ldots, m\\
        & \theta_k \in \cset_k,\quad k = 1, \ldots, K,
    \end{array}
\end{equation}
where the optimization variables are $\theta_k \in \reals^n$, $k = 1, \ldots, K$, and $z(t) \in \reals^K$, $t = 1, \ldots, m$.
Note that the variable $z(t)$ is a soft assignment vector for the hidden state label $\hat{z}(t)$, where the $k$th entry of $z(t)$ indicates the probability of the state being $k$ at the time step $t$, and $\hat{z}(t)$ can be estimated as the index of the largest entry of $z(t)$ after solving the problem (\ref{prob:iohmm}).
Each component of the problem (\ref{prob:iohmm}) can be interpreted as follows:
The first term in the objective function is the negative log-likelihood of the observed data under the IO-HMM model, given the state assignment probabilities $z(t)$, $t = 1, \ldots, m$, and the model parameters $\theta_k$, $k = 1, \ldots, K$.
The second term is a Tikhonov regularization on the model parameters $\theta_k$, with the regularization parameter $\alpha_\theta > 0$.
The third term is a temporal smoothness regularization on the state assignment probabilities, where $\Dkl(p, q)$ denotes the Kullback-Leibler divergence between two probability distributions $p$ and $q$, and $\alpha_z > 0$ is the corresponding regularization parameter.
The constraints on the variables $z(t)$, $t = 1, \ldots, m$, ensure that they are valid probability distributions.
The sets $\cset_k \subseteq \reals^n$, $k = 1, \ldots, K$, are closed nonempty convex sets that encode potential prior knowledge about the model parameters $\theta_k$.

In this example, we generate $m = 1800$ data samples from an IO-HMM with $K = 3$ hidden states and input feature dimension $n = 2$.
The feature vector for each sample is generated according to
\begin{equation*}
    x(t) \sim ({\cal U}(-5, 5),\ 1),
\end{equation*}
where ${\cal U}(a, b)$ denotes a uniform distribution over the interval $[a, b]$, and the second entry of $x(t)$ is always $1$ to account for the bias term.
The ground truth coefficients, initial state distribution, and transition matrix are given by 
\begin{align*}
    &\theta_1 = (-1, 0),\quad \theta_2 = (2, 6),\quad \theta_3 = (2, -6),\\
    &p_{\rm init} = (1, 0, 0),\quad 
    P_{\rm tr} = \left[\begin{array}{lll}
        0.95 & 0.025 & 0.025\\
        0.025 & 0.95 & 0.025\\
        0.025 & 0.025 & 0.95
    \end{array}\right].
\end{align*}
To fully specify the problem (\ref{prob:iohmm}), it is assumed that we are given the following prior knowledge about the coefficients:
\begin{equation*}
    \theta_{1,1} \leq 0,\quad \theta_{2, 1} \geq 0,\quad \theta_{3, 1} \geq 0,\quad \theta_{2, 2} \geq \theta_{3, 2},
\end{equation*}
where $\theta_{i, j}$ denotes the $j$th entry of the vector $\theta_i$.
The corresponding code snippet is as follows:
\begin{lstlisting}[language=zhpython]
thetas = cp.Variable((K, n))
zs = cp.Variable((m, K), nonneg=True)

rs = [
    -cp.multiply(ys, xs @ thetas[k])
    + cp.logistic(xs @ thetas[k])
    for k in range(K)
]
obj = cp.Minimize(
    cp.sum(cp.multiply(zs, cp.vstack(rs).T))
    + alpha_theta * cp.sum_squares(thetas)
    + alpha_z * cp.sum(cp.kl_div(zs[:-1], zs[1:])))
constr = [
    thetas[0][0] <= 0,
    thetas[1][0] >= 0,
    thetas[2][0] >= 0,
    thetas[1][1] >= thetas[2][1],
    zs <= 1, cp.sum(zs, axis=1) == 1
]

prob = BiconvexProblem(obj, [zs], [thetas], constr)
\end{lstlisting}
For this example, we set \texttt{mode='penalty'} when calling \texttt{BiconvexProblem.solve()}.
This solves the relaxed biconvex problem (\ref{prob:bcp_relax}) via the penalty ACS procedure (algorithm~\ref{alg:infeas_acs}).

\begin{figure}[t]
    \centering
    \includegraphics[width=0.8\textwidth]{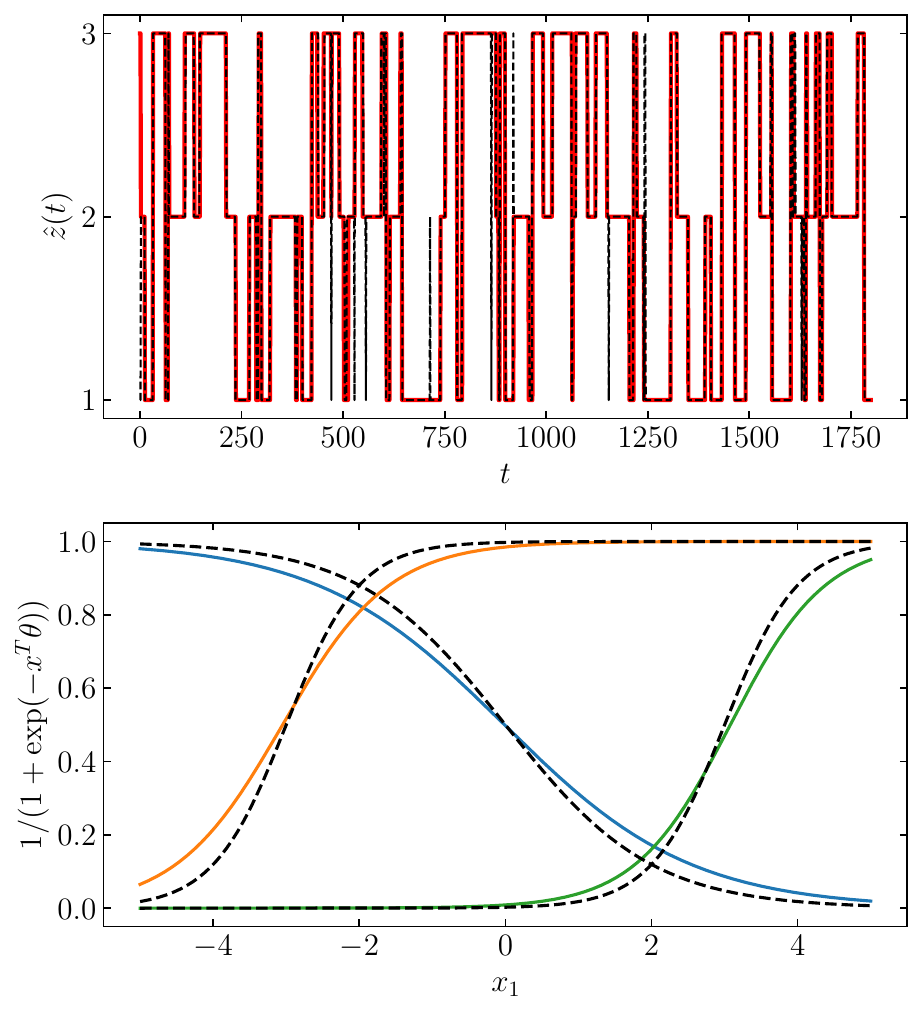}%
    \caption{
        Results of the IO-HMM example.
        The ground truth state labels (\emph{top}) and the response curve of the output probability (\emph{bottom}) are shown in black dashed lines, whereas the corresponding estimates are shown in solid lines.
    }\label{fig:iohmm}
\end{figure}

The problem (\ref{prob:iohmm}) is then solved with regularization parameters $\alpha_\theta = 0.1$ and $\alpha_z = 2$.
The arguments to the \texttt{prob.solve()} method are set to \texttt{mode='penalty'}, \texttt{nu=1e2}, \texttt{lbd=0.1}, and \texttt{abs\_tol=1e-3}.
The final total slack is around $4.45 \times 10^{-8}$, indicating that the returned point is feasible (within numerical roundoff error) for the original biconvex problem (\ref{prob:iohmm}).
The top panel of figure~\ref{fig:iohmm} shows the decoded state labels $\hat{z}(t) = \argmax_{k \in \{1, \ldots, K\}} z_k(t)$, $t = 1, \ldots, m$, obtained from the optimized soft assignments, with ties resolved by the smallest maximizing index. The bottom panel shows the fitted response curves.
The estimated transition matrix of the Markov chain is
\begin{equation*}
    \left[\begin{array}{lll}
        0.96 & 0.02 & 0.02\\
        0.02 & 0.96 & 0.02\\
        0.02 & 0.03 & 0.95
    \end{array}\right].
\end{equation*}

\section{Discussion}
\subsection{Conclusion}
In this paper, we introduce DBCP, a DSL for specifying and solving biconvex optimization problems.
When an optimization problem is compliant with the DBCP ruleset, its biconvex structure can be verified and the problem can be automatically split into two DCP-compliant convex subproblems, which are then solved using an ACS procedure.
We implement these ideas in the open-source Python package \texttt{dbcp}, an extension of CVXPY that allows users to express biconvex problems in a high-level, human-readable way while automatically handling the construction and solution of the convex subproblems.

We should note that DBCP cannot eliminate the solution difficulties intrinsic to biconvex optimization.
Its ruleset is sufficient but not necessary, so not every biconvex formulation is recognized.
For compliant problems, ACS does not generally guarantee global optimality, and its result can depend on initialization.
What DBCP can do is to reduce the expert barrier in problem modeling and implementation, and to provide a convenient interface for users to explore different biconvex formulations and solution procedures.

The introduced examples show that models from several application areas can be expressed in compact, readable code close to their mathematical formulations.
We hope this lowers the barrier to biconvex optimization in practice and encourages applications in new areas.

\subsection{Comparison with DMCP}\label{sec:comp_dmcp}
The DMCP framework~\cite{shen2017disciplined} can be considered an early prototype for the DBCP framework we present here.
In addition to being based on an early version of CVXPY and no longer maintained, it has a different modeling scope and implementation design.

First of all, DMCP is designed for multi-convex problems, while DBCP is designed for biconvex problems.
The focus on two variable groups allows DBCP to draw directly on the more developed theoretical foundation for biconvex programming, including the partial optimality properties and ACS convergence results reviewed in \S\ref{sec:bcp} and \S\ref{sec:bcp_solution}; see also Wendell and Hurter~Jr.~\cite{wendell1976minimization} and Gorski~\etal~\cite{gorski2007biconvex}.
Secondly, as far as we know, most multi-convex problems that appear in practical applications can be formulated as biconvex problems by grouping variables appropriately.
As a very basic example, the multi-convex problem
\begin{equation*}
    \begin{array}{ll}
        \mbox{minimize} & |z_1z_2 + z_2z_3 + z_3z_4|\\
        \mbox{subject to} & z_1 + z_2 + z_3 + z_4 = 1
    \end{array}
\end{equation*}
with variables $z = (z_1, z_2, z_3, z_4) \in \reals^4$ is multi-convex since, when any three variables are \emph{fixed}~\cite{shen2017disciplined}, the objective is the absolute value of an affine function of the remaining variable, and the constraint is affine.
However, this problem can be reformulated as a biconvex problem by grouping the variables as $x = (z_1, z_3) \in \reals^2$ and $y = (z_2, z_4) \in \reals^2$, which leads to the biconvex problem
\begin{equation*}
    \begin{array}{ll}
        \mbox{minimize} & \left|x^T \left[\begin{array}{cc}1 & 0\\1 & 1\end{array}\right] y\right|\\[15pt]
        \mbox{subject to} & \ones^T x + \ones^T y = 1.
    \end{array}
\end{equation*}
Such a variable transformation allows each ACS step to update several variables simultaneously.
When a compatible two-group representation is available, it is reasonable to expect that these joint updates may improve the efficiency and robustness of the solution procedure compared with treating each variable as a separate block as in the multi-convex case.

Regarding the implementation design, the DMCP framework registers a block coordinate descent \emph{solver} for a standard CVXPY problem.
In contrast, DBCP uses a specialized \texttt{BiconvexProblem} object.
This object-oriented design allows the modeled biconvex problem and its components to be represented explicitly and inspected directly by a user.
The two frameworks also handle initialization and solution procedures differently.
The DMCP framework uses user specified or random initial values and introduces penalized slack variables throughout its block updates.
By default, DBCP checks the initial point and, when necessary, attempts a separate feasibility search before applying ACS to the original, unrelaxed problem.
This separates feasibility recovery from objective optimization, while the slack-penalized ACS is used as an optional alternative.
We believe that our design choices fit better into the current CVXPY ecosystem, and we hope that they could provide the user with more flexibility and control over, \eg, the problem modeling, solution process, and results inspection.

\subsection{Comparison with DNLP}\label{sec:comp_dnlp}
In parallel with the development of this work, the DNLP~\cite{cederberg2026disciplined} framework was independently developed and released.
We view DBCP and DNLP as complementary modeling options rather than competing replacements.
Although their modeling scopes overlap, they expose different problem structures and use different local solution procedures.
For example, the nonnegative matrix factorization problem presented in \S\ref{sec:eg_nmf} can be specified using either framework.
Several applications in the DNLP examples, including circle packing, path planning, sparse signal recovery, and phase retrieval, also appear among the earlier disciplined convex-concave programming examples~\cite{shen2016disciplined}.
This overlap illustrates that the same application can be approached through different disciplined modeling frameworks, whose rules and solution procedures remain distinct.

DNLP permits general smooth nonlinear expressions together with structured compositions of nonsmooth convex and concave functions.
For a specified general optimization problem, it checks whether the problem satisfies the DNLP rules and, if so, canonicalizes it into a smooth nonlinear program and passes it to a supported nonlinear programming solver, such as IPOPT~\cite{wachter2006implementation}, KNITRO~\cite{byrd2006knitro}, COPT~\cite{ge2022cardinal}, and Uno~\cite{vanaret2026implementing,vanaret2026uno}.
DBCP instead applies ACS to its two DCP-compliant convex subproblems using CVXPY-compatible convex optimization solvers.
DBCP therefore does not introduce a dependency on a nonlinear programming solver.
For users who already have the standard CVXPY convex-solving stack installed, this results in a lighter set of solver dependencies.

The initialization procedures of the two frameworks also differ.
DNLP initializes auxiliary variables introduced during canonicalization and passes the complete nonlinear program, with user specified or generated initial values, to a nonlinear programming solver.
By comparison, the direct method in DBCP uses the separate feasibility search workflow described above before applying ACS to the original problem.
Since both frameworks use local solution procedures, they can be sensitive to initialization and may converge to different local optima.
In the most general case, because neither framework can guarantee global optimality, there is not much we can say about the relative quality of the returned solutions.

From the perspective of application, DNLP applies to a broader range of nonlinear formulations, while DBCP provides an alternative when a useful biconvex decomposition is available and the user wishes to expose and exploit that structure.
On the other hand, since both frameworks build on CVXPY variables, objectives, constraints, and expressions, when a formulation satisfies both rule systems, its objective and constraints can generally be reused by changing the problem construction and solve call.
This makes it straightforward for users to try both approaches on modeled problems in their shared scope.

\newpage
\section*{Acknowledgments}
This work has been funded as part of BrainLinks-BrainTools, which is funded by the Federal Ministry of Economics, Science and Arts of Baden-W\"urttemberg within the sustainability program for projects of the Excellence Initiative II, and CRC/TRR 384 ``IN-CODE''.
The authors thank the anonymous reviewers for their constructive comments and suggestions, which helped improve the manuscript as well as the associated software.

\newpage
\appendix
\section{Biconvex problem with generalized inequality constraints}\label{sec:gen_ineq_constr}
Consider a biconvex optimization problem with generalized inequality constraints of the form
\begin{equation}\label{prob:bcp_gen_ineq}
    \begin{array}{ll}
        \mbox{minimize} & f_0(x, y)\\
        \mbox{subject to} & f_i(x, y) \preceq_{\cone_i} 0,\quad i = 1, \ldots, m\\
        & h_i(x, y) = 0,\quad i = 1, \ldots, p,
    \end{array}
\end{equation}
where $f_0 \colon \xset \times \yset \to \reals$ is a biconvex objective, $h_i \colon \xset \times \yset \to \reals$, $i = 1, \ldots, p$, are biaffine equality constraint functions, and $f_i \colon \xset \times \yset \to \reals^{q_i}$, $i = 1, \ldots, m$, are biconvex inequality constraint functions with respect to proper cones $\cone_i \subseteq \reals^{q_i}$, $i = 1, \ldots, m$.
The ACS procedure (algorithm~\ref{alg:acs}) and the proximal regularizations given by (\ref{eq:prox_acs_x_update}) and (\ref{eq:prox_acs_y_update}) can be directly applied to the generalized inequality constrained biconvex problem (\ref{prob:bcp_gen_ineq}) without modification.
However, the initialization procedure (algorithm~\ref{alg:bcp_init}) and the penalty ACS procedure (algorithm~\ref{alg:infeas_acs}) need to be slightly modified to accommodate the generalized inequality constraints.
Specifically, the relaxed biconvex feasibility problem (\ref{prob:bcp_feas_relax}) is changed to
\begin{equation}\label{prob:bcp_feas_relax_gen_ineq}
    \begin{array}{ll}
        \mbox{minimize} & \ones^T s + \norm{t}_1\\
        \mbox{subject to} & s \succeq 0\\
        & f_i(x, y) \preceq_{\cone_i} s_i e_{\cone_i},\quad i = 1, \ldots, m\\
        & h_i(x, y) = t_i,\quad i = 1, \ldots, p,
    \end{array}
\end{equation}
where $e_{\cone_i} \succeq_{\cone_i} 0$ is any positive element of the proper cone $\cone_i$.
For example, for a second-order cone $\cone = \{(x, t) \in \reals^q \mid \norm{x}_2 \leq t\}$, we can choose $e_{\cone} = (0, 1) \in \reals^{q}$;
for a positive semidefinite cone $\cone = \symms_+^q$, we can choose $e_{\cone} = I \in \reals^{q \times q}$, where $I$ is the identity matrix.
Similarly, the relaxed biconvex optimization problem (\ref{prob:bcp_relax}) is changed to
\begin{equation}\label{prob:bcp_relax_gen_ineq}
    \begin{array}{ll}
        \mbox{minimize} & f_0(x, y) + \nu(\ones^T s + \norm{t}_1)\\
        \mbox{subject to} & s \succeq 0\\
        & f_i(x, y) \preceq_{\cone_i} s_i e_{\cone_i},\quad i = 1, \ldots, m\\
        & h_i(x, y) = t_i,\quad i = 1, \ldots, p.
    \end{array}
\end{equation}
We implement the \texttt{dbcp} package in a way that these generalized inequality constraints are natively supported.
When a user specifies a DBCP-compliant biconvex problem, \texttt{dbcp} will automatically detect the presence of generalized inequality constraints, then determine whether the problem formulations (\ref{prob:bcp_feas_relax_gen_ineq}) and (\ref{prob:bcp_relax_gen_ineq}) should be used instead of (\ref{prob:bcp_feas_relax}) and (\ref{prob:bcp_relax}), and finally, modify the initialization procedure (algorithm~\ref{alg:bcp_init}) and the penalty ACS procedure (algorithm~\ref{alg:infeas_acs}) accordingly.
From a user's perspective, all these modifications are made in the backend and no extra input is required.

\newpage
\bibliography{refs}

\end{document}